\documentclass[letterpaper, 10 pt, conference]{ieeeconf}  

\IEEEoverridecommandlockouts                              

\usepackage[utf8]{inputenc}
\usepackage{float,color,cite}
\usepackage{dsfont}
\usepackage{bm}
\usepackage{stmaryrd}
\usepackage{booktabs,tabularx, hhline,multicol}
\usepackage{amsmath}
\usepackage{amssymb,amsfonts}
\usepackage{textcomp, gensymb}
\usepackage{array}
\usepackage[ruled,vlined,linesnumbered]{algorithm2e}
\usepackage{algorithmic}
\usepackage{hyperref,xcolor}
\SetKwInOut{Parameter}{Parameters}
\usepackage{graphicx}
\usepackage{booktabs}
\usepackage{diagbox}
\usepackage{autobreak}
\usepackage{comment}
\usepackage[english]{babel}
\newtheorem{theorem}{Theorem}

\newcommand{\enumbracket}[1]{{\llbracket #1 \rrbracket}}

\newcommand{\Rmnum}[1]{\uppercase\expandafter{\romannumeral #1}}

\ifodd 1

    \newcommand{\com}[2]{\textbf{\color{blue} (COMMENT from [#1]: #2)}}
\else
    \newcommand{\rev[1]}{{#1}}
    
    \newcommand{\com}[2]{}
\fi

\newcommand{\modelname}{PredVAR}

\makeatletter
\newcommand*\bigcdot{\mathpalette\bigcdot@{.5}}
\newcommand*\bigcdot@[2]{\mathbin{\vcenter{\hbox{\scalebox{#2}{$\m@th#1\bullet$}}}}}
\makeatother

\title{Probabilistic Reduced-Dimensional Vector Autoregressive Modeling for  Dynamics Prediction and Reconstruction with Oblique Projections }

\author{Yanfang Mo,~\IEEEmembership{Member,~IEEE}, Jiaxin Yu,~\IEEEmembership{Student Member, IEEE}, and~S.~Joe~Qin,~\IEEEmembership{Fellow,~IEEE}
\thanks{
Yanfang Mo is with the Hong Kong Institute for Data Science, City University of Hong Kong, Hong Kong {\tt\small Yanfang.MO@cityu.edu.hk}}
\thanks{
Jiaxin Yu is with the School of Data Science and Hong Kong Institute for Data Science,
City University of Hong  Kong, Hong Kong.
}
\thanks{S. Joe Qin is with the Department of Computing and Decision Science and Institute of Data Science, Lingnan University, Hong Kong. Corresponding author.
        {\tt\small joeqin@LN.edu.hk}. He acknowledges the financial support from  a Natural Science Foundation of China Project (U20A20189), a General Research Fund by RGC of Hong Kong (No. 11303421),
a Collaborative Research Fund by RGC of Hong Kong (Project No. C1143-20G),
        an ITF - Guangdong-Hong Kong Technology Cooperation Funding Scheme  (Project Ref. No. GHP/145/20), an InnoHK initiative of The Government of the
HKSAR for the Laboratory for AI-Powered Financial Technologies, and  a Shenzhen-Hong Kong-Macau Science and Technology Project Category C (9240086).}
}
\begin{document}


\maketitle
\thispagestyle{empty}
\pagestyle{empty}

\begin{abstract}
    In this paper, we propose a probabilistic reduced-dimensional vector autoregressive (PredVAR) model with oblique projections. This model partitions the measurement space into a dynamic  subspace and a static  subspace that do not need to be orthogonal. The partition allows us to apply an oblique projection to extract dynamic latent variables (DLVs) from high-dimensional data with maximized predictability. We develop an alternating iterative \modelname~algorithm that exploits the interaction between updating the latent VAR dynamics and estimating the oblique projection, using expectation maximization (EM) and a statistical constraint. In addition, the noise covariance matrices are estimated as a natural outcome of the EM method. A simulation case study  of the nonlinear Lorenz oscillation system illustrates the advantages of the proposed approach over two alternatives.
\end{abstract}

\section{Introduction}
Extracting reduced-dimensional dynamics is crucial in many industries, including chemicals, power systems, finance, and transportation, where operational data are usually high-dimensional with dynamic features~\cite{sznaier2020control,pena1987identifying,qin2020bridging,gao2021modeling,reinsel2023multivariate}.
Classic data analytic tools like principal component analysis (PCA) and canonical correlation analysis (CCA) are oblivious to the system dynamics~\cite{qin2020bridging} and should be extended. Dynamic PCA (DPCA)~\cite{ku1995disturbance} and subspace-based models~\cite{li2001consistent} can consider the auto-correlations in time series. A new dynamic latent variable (DLV) model was developed in~\cite{li2014new} to characterize the dynamic relations in DLVs and the static cross-correlations in residuals. Moreover, dynamic-inner PCA (DiPCA)~\cite{dong2018novel} and dynamic-inner CCA (DiCCA)~\cite{dong2018dynamic} algorithms produce rank-ordered DLVs by maximizing prediction power. Furthermore, Qin developed a latent vector autoregressive modeling algorithm with a CCA objective (LaVAR-CCA)~\cite{qin2021latent} and a vector autoregressive model~\cite{qin2022latent} for the DLVs, where a state-space generalization is subsequently developed in~\cite{yu2022latent}.

Uncertainty estimation requires a statistical viewpoint for dimension reduction and dynamics extraction. The LaVAR-CCA~\cite{qin2022latent} has a statistical interpretation of profile likelihood. Many dynamic factor models (DFMs) use time-related statistics to estimate model parameters~\cite{box1977canonical,pena1987identifying,lam2012factor,pena2019forecasting,gao2021modeling}. However, more attention should be paid to dynamics modeling in the statistical literature~\cite{qin2022latent,wen2012data,zhou2016autoregressive}. More importantly, the structured dynamics-noise relationship can facilitate signal reconstruction and prediction~\cite{gao2021modeling,qin2020bridging}.

Our probabilistic model partitions the measurement space into a signal subspace admitting the low-dimensional DLV dynamics and a static noise subspace, which do not need to be orthogonal to each other. Thus, an oblique-projection perspective is adopted for the signal-noise structure. This perspective is well-known in the signal-processing literature~\cite{behrens1994signal}, where model parameters are deemed known, which is not the case for this work. A challenge lies in identifying a proper oblique projection to extract DLVs and degrade the effect of noise on estimating the system dynamics.

Also, the estimated DLV dynamics can facilitate oblique projection identification. Thus, an alternating iterative scheme is developed for the latent dynamics and oblique projection estimations by exploring their interactions. Specifically, two expectation-maximization (EM) steps are used to estimate the DLV dynamics and the signal subspace. Meanwhile, a statistical constraint is uniquely imposed to characterize the relationship between the estimated signal and static noise and is thus used to identify the noise subspace.  The contributions in this work are as follows.

$1)$ A probabilistic reduced-dimensional vector regressive (PredVAR) model with oblique projections is proposed. 

$2)$ An iterative algorithm is developed to update the DLV dynamics and oblique projection estimations alternately. It uniquely uses a statistical constraint together with an EM procedure in the oblique projection identification.

$3)$ A simulated case study is conducted to demonstrate the strength of our approach compared to two benchmarks, a one-shot algorithm that first identifies the oblique projection and then estimates the DLV dynamics and a counterpart that focuses on the orthogonal projection.



\section{Model Formulation}\label{Sec_Model}
\begin{figure}[t]
    \centering
\includegraphics[width=0.7\columnwidth]{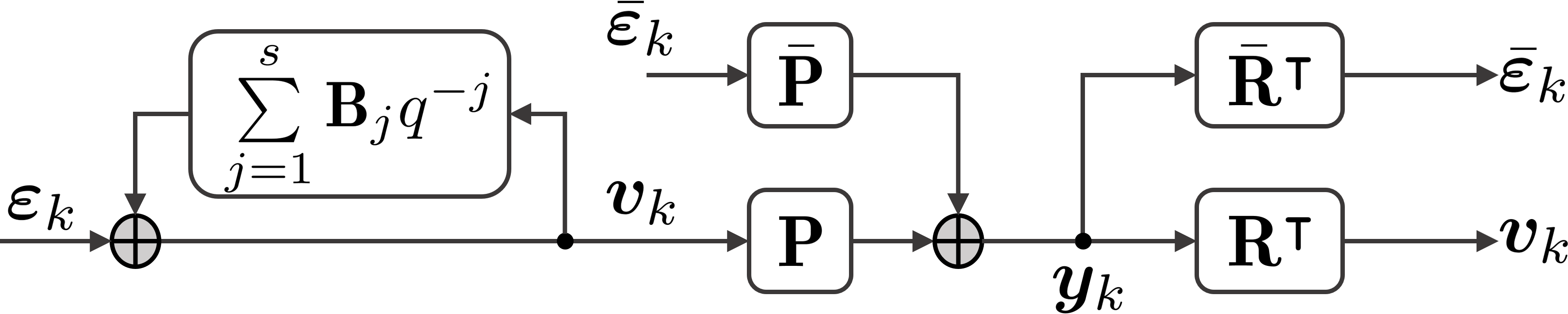}
    \caption{A block diagram of a \modelname~model.}
    \label{fig:block_diag}
\end{figure}

Denote the measurement time series from a dynamic system by~$\{\bm y_k \in \Re^p\}_{k=1}^{N+s}$. It is typical for operational data that the dynamics are excited in a subspace of dimension~$\ell < p$, referred to as reduced-dimensional dynamics~\cite{qin2022latent}. Thus, a dynamic latent process~$\{\bm{v}_k\in \Re^\ell \}$ exists, together with a loadings matrix~$\mathbf{P} \in \Re^{p\times \ell }$ that is of full column rank.

The reduced-dimensional dynamics and the signal-noise structure motivate us to construct the model in~\cite{qin2022latent} in the \modelname~model as
\begin{align}
    &\bm y_k = \mathbf{P} \bm v_k + \mathbf{\bar{P}} \bm{\bar \varepsilon}_k,  &\bm{\bar \varepsilon}_k \sim \mathcal{N}(\bm 0, \mathbf{\Sigma}_{\bm{\bar \varepsilon}}),\label{eq:outer_model}\\
    &\bm v_k = \sum_{j=1}^s\mathbf{B}_j\bm v_{k-j} + \bm{\varepsilon}_k,~~ &\bm{\varepsilon}_k \sim \mathcal{N}(\bm 0, \mathbf{\Sigma}_{\bm \varepsilon}).\label{model_dy}
\end{align}
where $[\mathbf{P} \quad \mathbf{\bar P}]$ is required to be nonsingular, the innovations vector $\bm \varepsilon_k\in \Re^{\ell}$ and the static noise~$\bm{\bar \varepsilon}_k \in \Re^{p-\ell}$ are assumed to be serially and mutually independent ($E\{\bm{\varepsilon}_k\bm{\bar \varepsilon}_k^\intercal\}=\bm 0$). An~$s$-order VAR model is used to describe the DLV dynamics. Equations \eqref{model_dy} and \eqref{eq:outer_model} are called inner and outer models. By~\eqref{model_dy}, the one-step-ahead prediction of~${\bm v}_k$ is defined as~
\begin{multline} \label{one-step-pred}
    \bm{\tilde v}_k =E\{\bm v_k\mid \bm y_1, \ldots, \bm y_{k-1}\}
     = \sum_{j=1}^s\mathbf{B}_j\bm v_{k-j}.
\end{multline}

Associated with~$\mathbf{P}$ and~$\mathbf{\bar{P}}$, the DLV and static weight matrices~$\mathbf{R}\in \Re^{p\times \ell }$ and~$\mathbf{\bar{R}}\in \Re^{p\times (p-\ell )}$ satisfy
\begin{equation} \label{eq:RP_relation}
    \begin{bmatrix}
        \mathbf{R} & \mathbf{\bar R}
    \end{bmatrix}^\intercal
    \begin{bmatrix}
        \mathbf{P} & \mathbf{\bar P}
    \end{bmatrix}
    = \mathbf{I}.
\end{equation}
Therefore, pre-multiplying $\mathbf{R}^\intercal$ to \eqref{eq:outer_model} and~$\mathbf{P}$ to~\eqref{eq:vk} lead to
\begin{equation} \label{eq:vk}
 \bm v_k = \mathbf{R}^\intercal\bm y_k \text{ and } \mathbf{P}\bm v_k = \mathbf{P}\mathbf{R}^\intercal \bm y_k,
\end{equation}
giving an oblique projection of $\bm y_k$ since $(\mathbf{P}\mathbf{R}^\intercal)^2 = \mathbf{P}\mathbf{R}^\intercal$.

The \modelname~model is illustrated in Fig.~\ref{fig:block_diag}. The VAR model for DLVs describes the low-dimensional dynamics and initiates exploring more complicated dynamics descriptions (see~\cite{dong2019efficient, yu2022latent,pillonetto2014kernel,khosravi2023existence}). By~\eqref{eq:outer_model}, \eqref{model_dy}, and \eqref{eq:vk}, a VAR model with reduced-rank coefficients~$\{\mathbf{P}\mathbf{B}_j\mathbf{R}^\intercal,j\in\enumbracket{s} \triangleq\{1,2,\ldots, s\}\}$ is obtained
\begin{align}
& \bm y_k =\sum_{j=1}^s\mathbf{P}\mathbf{B}_j\mathbf{R}^\intercal \bm y_{k-j} +\bm e_k\text{, where } \label{model_measurement}\\
 &  \bm e_k = \mathbf{P}\bm \varepsilon_k +\mathbf{\bar{P}}\bm{\bar \varepsilon}_k \sim \mathcal{N}(\bm 0, \mathbf{\Sigma}_{\bm e}=\mathbf{P}\mathbf{\Sigma}_{\bm\varepsilon}\mathbf{P}^{\intercal}+\mathbf{\bar{P}}\mathbf{\Sigma}_{\bm{\bar \varepsilon}}\mathbf{\bar{P}}^{\intercal}).\label{eq:VAR_yk_noise}
\end{align}

By~\eqref{eq:RP_relation}, \eqref{eq:VAR_yk_noise}, and the assumption~$E\{\bm{\varepsilon}_k\bm{\bar \varepsilon}_k^\intercal\}=\bm 0$, the  constraint below is derived and imposed in this paper
\begin{equation}\label{constraint:PR}
\mathbf{R}^\intercal\mathbf{\Sigma}_{\bm e}\mathbf{\bar{R}} = \mathbf{0}.
\end{equation}

The imposed statistical constraint~\eqref{constraint:PR} plays a role in estimating the parameter tuple~$(\mathbf{P};\bm v;\mathbf{\Sigma}_{\bm \varepsilon};  \mathbf{B}_j,j\in \enumbracket{s};\mathbf{\bar{P}};\mathbf{\Sigma}_{\bm{\bar \varepsilon}})$, discussed in the next section. Two tuples of \modelname~model parameters are observationally equivalent if they can generate the same measurement data. In this regard, the following theorem deserves attention when estimating a \modelname~model with given measurements~\cite{bai2012statistical}. 

\begin{theorem}\label{thm:eq_two_models}
For arbitrary nonsingular matrices~\mbox{$\mathbf{M}\in\Re^{\ell \times \ell }$} and~$\mathbf{\bar{M}}\in\Re^{(p-\ell )\times (p-\ell )}$, the two \modelname~models are observationally equivalent with their respective parameter tuples as~$(\mathbf{P};\bm v;\mathbf{\Sigma}_{\bm \varepsilon};  \mathbf{B}_j,j\in \enumbracket{s};\mathbf{\bar{P}};\mathbf{\Sigma}_{\bm{\bar \varepsilon}})$ and~$(\mathbf{P}\mathbf{M}^{-1};\mathbf{M}\bm v;$ $\mathbf{M}\mathbf{\Sigma}_{\bm \varepsilon}\mathbf{M}^{\intercal};\mathbf{M}\mathbf{B}_j\mathbf{M}^{-1},j\in \enumbracket{s};\mathbf{\bar{P}}\mathbf{\bar{M}}^{-1};\mathbf{\bar{M}}\Sigma_{\bm{\bar \varepsilon}}\mathbf{\bar{M}}^{\intercal}).$
\end{theorem}

Note that the oblique projections related to~$\mathbf{P}$ and~$\mathbf{R}$ are the same for observationally equivalent models in Theorem~\ref{thm:eq_two_models}~\cite{gao2021modeling}. So are the matrices~$\mathbf{P}\mathbf{\Sigma}_{\bm \varepsilon}\mathbf{P}^{\intercal}$,~$\mathbf{\bar{P}}\mathbf{\Sigma}_{\bm{\bar \varepsilon}}\mathbf{\bar{P}}^{\intercal}$,~$\mathbf{\Sigma}_{\bm e}$, and~$\mathbf{P}\mathbf{B}_j\mathbf{R}^\intercal,j\in\enumbracket{s}$. Moreover, diverse restrictions can be imposed to improve the identifiability and for numerical concerns, as shown in~\cite{bai2012statistical}.

\section{Model Estimation }\label{Sec_Alg}
\begin{figure}[t]
    \centering    \includegraphics[width=0.62\columnwidth]{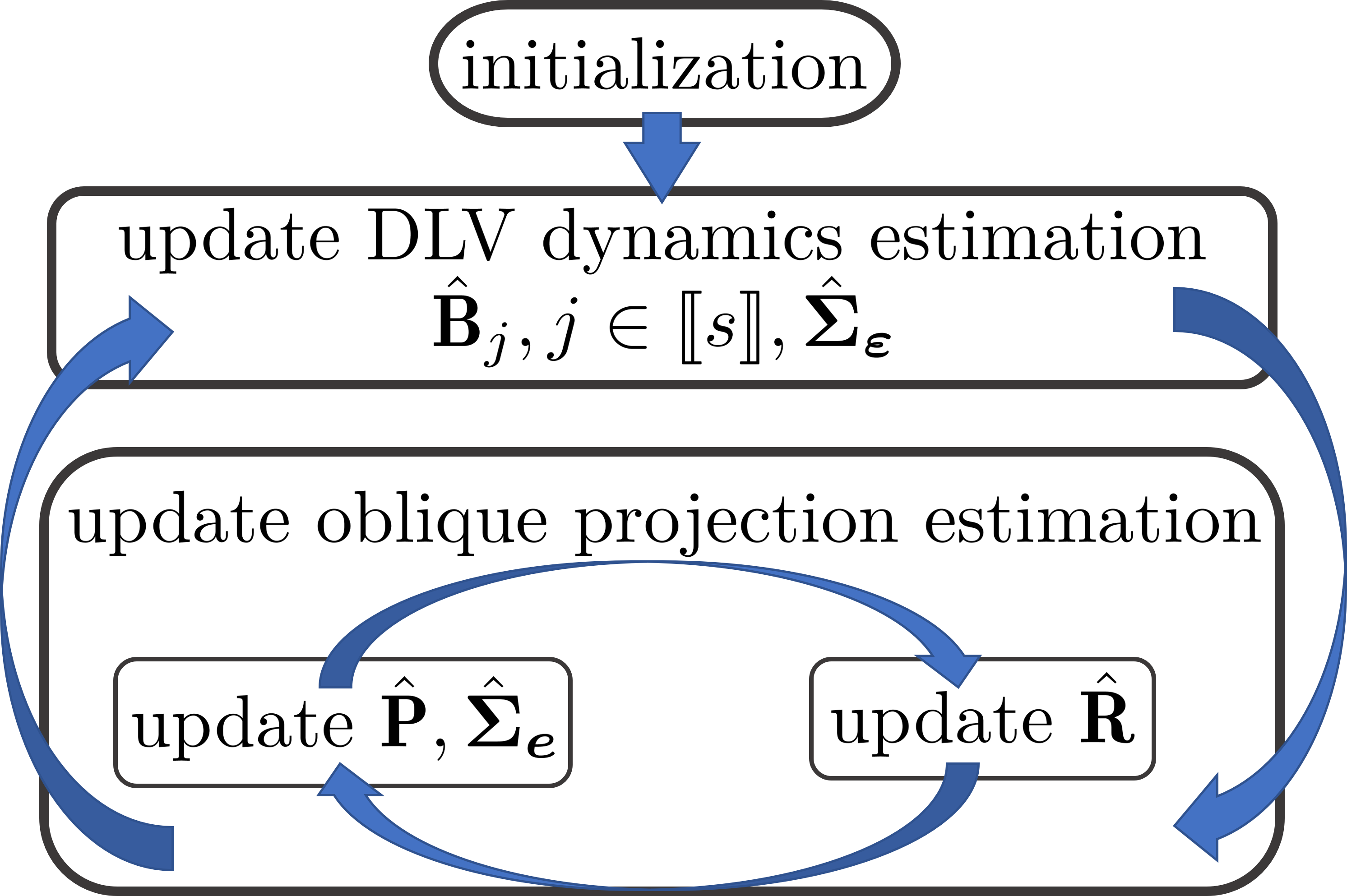}
    \vspace{-0.5em}
    \caption{An illustration of our alternating procedure.}
    \label{fig:flow_chart}
    \vspace{-1.6em}
\end{figure}

To identify the \modelname~model, the DLV dynamics and the oblique projection estimations are alternately conducted via EM and a statistical constraint, as depicted in Fig.~\ref{fig:flow_chart}. For convenience, the symbol~$\hat{\bigcdot}$ is used to distinguish between an estimated parameter and its true value.

\subsection{Identifying DLV Dynamics with Oblique Projections}\label{Sec._Id_DLV}

Given~$(\mathbf{P},\mathbf{\bar{P}})$, to use EM for estimating~$\mathbf{\Sigma}_{\bm \varepsilon}$ and~$\mathbf{B}_j,j\in \enumbracket{s}$, the following likelihood should be maximized:
\begin{equation*}
     \prod_{k=s+1}^{s+N} p(\bm v_{k}\mid \bm v_{k-1}, \bm v_{k-2}, \ldots,\bm v_{k-s}).
\end{equation*}
Since~$\bm{\varepsilon}_k \sim \mathcal{N}(\bm 0,\mathbf{\Sigma}_{\bm \varepsilon})$, the M-step requires minimizing the following function
\begin{multline*}
    L^{DLV}(\mathbf{\Sigma}_{\bm \varepsilon},\mathbb{B}) =N\ln|\mathbf{\Sigma}_{\bm \varepsilon}| +\\\sum_{k=s+1}^{s+N}(\bm v_{k}-\bm{\tilde v}_{k})^\intercal \mathbf{\Sigma}_{\bm \varepsilon}^{-1} (\bm v_{k}-\bm{\tilde v}_{k}),
\end{multline*}
where~$\bm{\tilde v}_{k}$ is the one-step-ahead prediction defined in~\eqref{one-step-pred}.
Take the derivatives of~$L^{DLV}(\mathbf{\Sigma}_{\bm \varepsilon},\mathbb{B})$ with respect to~$\mathbf{\Sigma}_{\bm \varepsilon}^{-1}$ and~$\mathbf{B}_j, j\in\enumbracket{s}$ and set them to zeros. It follows that
\begin{align*}
    &\mathbf{\Sigma}_{\bm \varepsilon} = \frac{1}{N}\sum_{k=s+1}^{s+N}(\bm v_{k}-\bm{\tilde v}_{k})(\bm v_{k}-\bm{\tilde v}_{k})^\intercal,\\
    &\mathbf{B}_j= \sum_{k=s+1}^{s+N}\Bigg(\bm v_{k}\bm v_{k-j}^\intercal-\!\!\sum_{i\in\enumbracket{s},i\neq j}\!\mathbf{B}_i\bm v_{k-i}\bm v_{k-j}^\intercal\Bigg)\times\\
    &\hspace{12em}
    \Bigg(\sum_{k=s+1}^{s+N} \bm v_{k-j} \bm v_{k-j}^\intercal \Bigg)^{-1}.
\end{align*}
 The last equation is rearranged as
\[
\sum_{i\in\enumbracket{s}} \mathbf{B}_i
\sum_{k=s+1}^{s+N} \bm v_{k-i}\bm v_{k-j}^\intercal= \sum_{k=s+1}^{s+N} \bm v_{k} \bm v_{k-j}^\intercal, \quad j\in\enumbracket{s}.
\]

In the E-step, by~\eqref{eq:vk}, the following updating formula is obtained with an estimate~$\hat{\mathbf{R}}$: \newline
 \smallskip \centerline{\medskip $\hat{\bm v}_k = \hat{\mathbf{R}}^{\intercal} \bm y_k.$}

Form the following augmented matrices as in~\cite{qin2022latent}:
\begin{align*}
    &\mathbf{Y}_i = [\bm y_{i+1}~\bm y_{i+2}~\cdots~\bm y_{i+N}]^{\intercal}, i\in \{0\}\cup\enumbracket{s};\\
    &\hat{\mathbf{V}}_i = [\hat{\bm v}_{i+1}~\hat{\bm v}_{i+2}~\cdots~\hat{\bm v}_{i+N}]^{\intercal}, i\in \{0\}\cup\enumbracket{s};\\
    & \hat{\mathbb{V}}= [\hat{\mathbf{V}}_{s-1}~\hat{\mathbf{V}}_{s-2}~\cdots~\hat{\mathbf{V}}_0];\\
    & \hat{\mathbb{B}} = [\hat{\mathbf{B}}_1~\hat{\mathbf{B}}_2~\cdots~\hat{\mathbf{B}}_s]^{\intercal}.
\end{align*}
Then, the previous formulas can be rewritten as
\begin{align}
&\hat{\mathbf{V}}_i=\mathbf{Y}_i\hat{\mathbf{R}}, i\in\{0\}\cup\enumbracket{s} \label{update_DLV}\\
    &\hat{\mathbb{B}} = (\hat{\mathbb{V}}^{\intercal}\hat{\mathbb{V}})^{-1}\hat{\mathbb{V}}^{\intercal}\hat{\mathbf{V}}_s; \label{update_B}\\
    & \hat{\mathbf{\Sigma}}_{\bm \varepsilon} = (\hat{\mathbf{V}}_s-\hat{\mathbb{V}}\hat{\mathbb{B}})^{\intercal}(\hat{\mathbf{V}}_s-\hat{\mathbb{V}}\hat{\mathbb{B}})/N \label{update_dn}.
\end{align}

The updating formula for~$\mathbb{B}$ is the same as that in LaVAR-CCA. This is unsurprising because LaVAR-CCA has a likelihood interpretation as stated in~\cite{qin2022latent}.

\subsection{Identifying Oblique Projection with DLV Dynamics}\label{Sec._Id_OB}
Given the DLV dynamics, to use EM for estimating~$\hat{\mathbf{P}}$ and~$\hat{\mathbf{R}}$, the following likelihood is maximized:
\begin{equation*}
     \prod_{k=s+1}^{s+N} p(\bm y_{k}\mid \bm v_{k-1}, \bm v_{k-2}, \ldots, \bm v_{k-s}),
\end{equation*}
while respecting~\eqref{eq:VAR_yk_noise} and~\eqref{constraint:PR} for~$\bm e_k$. This is equivalent to minimizing the function below, subject to~\eqref{constraint:PR}.
\begin{multline*}
    L^{proj}(\mathbf{\Sigma}_{\bm e},\mathbf{P}) =N\ln|\mathbf{\Sigma}_{\bm e}| +\\\sum_{k=s+1}^{s+N}(\bm y_{k}-\mathbf{P}\bm{\tilde v}_{k})^\intercal \mathbf{\Sigma}_{\bm e}^{-1} (\bm y_{k}-\mathbf{P}\bm{\tilde v}_{k}),
\end{multline*}


The imposed constraint~\eqref{constraint:PR} allows us to uniquely determine the subspace spanned by~$\mathbf{R}$, given~$\mathbf{P}$ and~$\mathbf{\Sigma}_{\bm e}$. Precisely, with singular value decompositions (SVDs), obtain~$\mathbf{\bar{R}}$ from an orthonormal basis of the null space of~$\mathbf{P}$ and then an orthonormal basis of the null space of~$\mathbf{\Sigma}_{\bm e}\mathbf{\bar{R}}$ gives~$\mathbf{R}$. 

Ignoring~\eqref{constraint:PR}, an EM procedure can be used to estimate~$\mathbf{P}$ and~$\mathbf{\Sigma}_{\bm e}$. Similarly to before, take the derivatives of~$L^{proj}(\mathbf{\Sigma}_{\bm e},\mathbf{P})$ with respect to $\mathbf{P}$ and $\mathbf{\Sigma}_{\bm e}^{-1}$  and set them to zero. Then, the M-step requires the updating formulas
\begin{align} &\hat{\mathbf{P}}=\mathbf{Y}_s^{\intercal}\hat{\mathbb{V}}\hat{\mathbb{B}}(\hat{\mathbb{B}}^{\intercal}\hat{\mathbb{V}}^{\intercal}\hat{\mathbb{V}}\hat{\mathbb{B}})^{-1}; \label{update_P}\\
    & \hat{\mathbf{\Sigma}}_{\bm e} = (\mathbf{Y}_s-\hat{\mathbb{V}}\hat{\mathbb{B}}\hat{\mathbf{P}}^{\intercal})^{\intercal}(\mathbf{Y}_s-\hat{\mathbb{V}}\hat{\mathbb{B}}\hat{\mathbf{P}}^{\intercal})/N \label{update_n}.
\end{align}
The E-step requires~\eqref{update_DLV}. The above analysis leads to Algorithm~\ref{alg: OB}, where~$\hat{\mathbf{P}}$ and~$\hat{\mathbf{R}}$ are alternately updated by the EM method and the statistical constraint~\eqref{constraint:PR} until convergence. Note that \eqref{update_P} is different from that in LaVAR-CCA~\cite{qin2022latent}, which updates~$\hat{\mathbf{P}}$ by the formula~$\hat{\mathbf{P}}=\mathbf{Y}_s^{\intercal}\hat{\mathbf{V}}_s(\hat{\mathbf{V}}_s^{\intercal}\hat{\mathbf{V}}_s)^{-1}$ and does not directly involve the DLV dynamics.

\begin{algorithm}[t]
\caption{Update~$\hat{\mathbf{P}}$ and~$\hat{\mathbf{R}}$ with~$\hat{\mathbf{B}}$}\label{alg: OB}
\LinesNumbered
\SetNoFillComment
\KwIn{$\{\bm y_k \in \Re^p\}_{k=1}^{N+s}$; $\hat{\mathbf{B}}_j,j\in \enumbracket{s}$; $\hat{\mathbf{R}}$;}
  \KwOut{$\hat{\mathbf{\Sigma}}_{\bm e}$;  $\hat{\mathbf{P}}$; $\hat{\mathbf{R}}$ with~$\hat{\mathbf{R}}^\intercal\hat{\mathbf{R}}=\mathbf{I}$; }
  \While {the termination condition is unsatisfied}{
   {Update~$\hat{\mathbf{P}}$ in~\eqref{update_P} and~$\hat {\mathbf{\Sigma}}_{\bm e}$ in~\eqref{update_n}\;}
   {Perform SVD on~$\hat{\mathbf{P}}=\mathcal{U}_1\mathcal{D}_1\mathcal{V}_1^\intercal$ with singular} {values in ascending order; ${\mathbf{\hat{\bar R}}}=\mathcal{U}_1(:,1:(p-\ell ))$\;}
   {Perform SVD on~$\hat{\mathbf{\Sigma}}_{\bm e}{\mathbf{\hat{\bar R}}}=\mathcal{U}_2\mathcal{D}_2\mathcal{V}_2^\intercal$ with singular} {values in ascending order;~$\hat{\mathbf{R}}=\mathcal{U}_2(:,1:\ell )$\;\label{alg: updateR}}
  {Calculate the DLVs in~\eqref{update_DLV}\;}
  }
\end{algorithm}

\begin{algorithm}[t]
\caption{A \modelname~Algorithm}\label{alg: DP}
\LinesNumbered
\SetNoFillComment
\KwIn{scaled measurements $\{\bm y_k \in \Re^p\}_{k=1}^{N+s}$;~$s$;~$l$; }
  \KwOut{$\hat{\mathbf{\Sigma}}_{\bm \varepsilon}$;  $\hat{\mathbf{B}}_j,j\in \enumbracket{s}$; $\hat{\mathbf{R}}$;~$\hat{\mathbf{P}}$;~$\hat{\mathbf{\Sigma}}_{\bm e}$;}
  {Perform SVD on~${\mathbf{Y}}_s^{\intercal}{\mathbf{Y}}_s/N=\mathcal{U}_0\mathcal{D}_0\mathcal{V}_0^\intercal$; }{$\hat{\mathbf{{R}}}=\mathcal{U}_0(:,1:({p-\ell }))$\;}
  \While {the termination condition is unsatisfied}{
  {Estimate the DLVs in~\eqref{update_DLV}\; \label{alg: updateDLV}}
  {Update~$\mathbb{B}$ in~\eqref{update_B} and~$\mathbf{\Sigma}_{\bm \varepsilon}$ in~\eqref{update_dn}\; \label{alg: updateB}}
    {Update~$\hat{\mathbf{R}}$,~$\hat{\mathbf{P}}$, and~$\hat{\mathbf{\Sigma}}_{\bm e}$ in Algorithm~\ref{alg: OB}\;}
  }
\end{algorithm}


\subsection{A Dynamics-Projection Alternate Updating Scheme}
Exploiting the interaction between the latent dynamics estimation in Section~\ref{Sec._Id_DLV} and the oblique projection estimation in Section~\ref{Sec._Id_OB}, we develop an alternating iterative scheme to identify a \modelname~model with an illustration in Fig.~\ref{fig:flow_chart} and pseudo-codes in Algorithm~\ref{alg: DP}. The algorithm's convergence follows from the classic result on EM~\cite{wu1983convergence} if the~$\mathbf{R}$-estimates are replaced with true values. Moreover, the estimate~$\hat{\mathbf{R}}$ is imposed to be orthonormal and solely depends on the statistical constraint~\eqref{constraint:PR} given other parameters.

Note that, although~$\mathbf{R}^\intercal\mathbf{P}=\mathbf{I}$ is not enforced in the algorithm, it follows from~\eqref{update_DLV},~\eqref{update_B}, and~\eqref{update_P} that
\begin{align*}
    \hat{\mathbf{R}}^\intercal \hat{\mathbf{P}}&\stackrel{\eqref{update_P}}{=}\hat{\mathbf{R}}^\intercal\mathbf{Y}_s^{\intercal}\hat{\mathbb{V}}\hat{\mathbb{B}}(\hat{\mathbb{B}}^{\intercal}\hat{\mathbb{V}}^{\intercal}\hat{\mathbb{V}}\hat{\mathbb{B}})^{-1}\stackrel{\eqref{update_DLV}}{=}\mathbf{V}_s^{\intercal}\hat{\mathbb{V}}\hat{\mathbb{B}}(\hat{\mathbb{B}}^{\intercal}\hat{\mathbb{V}}^{\intercal}\hat{\mathbb{V}}\hat{\mathbb{B}})^{-1}\\
    &\stackrel{\eqref{update_B}}{=}\mathbf{V}_s^{\intercal}\hat{\mathbb{V}}(\hat{\mathbb{V}}^\intercal\hat{\mathbb{V}})^{-1}\mathbf{V}_s(\mathbf{V}_s^{\intercal}\hat{\mathbb{V}}(\hat{\mathbb{V}}^\intercal\hat{\mathbb{V}})^{-1}\mathbf{V}_s)^{-1}=\mathbf{I},
\end{align*}
and follows from~\eqref{update_DLV},~\eqref{update_dn},~\eqref{update_n}, and~$\hat{\mathbf{P}}^\intercal \hat{\mathbf{R}}=\mathbf{I}$ that
\begin{align*}
    \hat{\mathbf{R}}^\intercal \hat{\mathbf{\Sigma}}_{\bm e} \hat{\mathbf{R}} &\stackrel{\eqref{update_n}}{=}  \hat{\mathbf{R}}^\intercal ({\mathbf{Y}}_s-\hat{\mathbb{V}}\hat{\mathbb{B}}\hat{\mathbf{P}}^{\intercal})^{\intercal}({\mathbf{Y}}_s-\hat{\mathbb{V}}\hat{\mathbb{B}}\hat{\mathbf{P}}^{\intercal})\hat{\mathbf{R}}/N\\
    &\stackrel{\eqref{update_DLV}}{=}(\hat{\mathbf{V}}_s-\hat{\mathbb{V}}\hat{\mathbb{B}})^{\intercal}(\hat{\mathbf{V}}_s-\hat{\mathbb{V}}\hat{\mathbb{B}})/N\stackrel{\eqref{update_dn}}{=} \hat{\mathbf{\Sigma}}_{\bm \varepsilon}
\end{align*}
after the convergence of Algorithm~\ref{alg: DP}.

\section{Simulation Case Study} \label{Sec_Simu}
\begin{figure*}[t]
\centering
\begin{minipage}{.5\textwidth}
    \centering
    \includegraphics[scale=0.40]{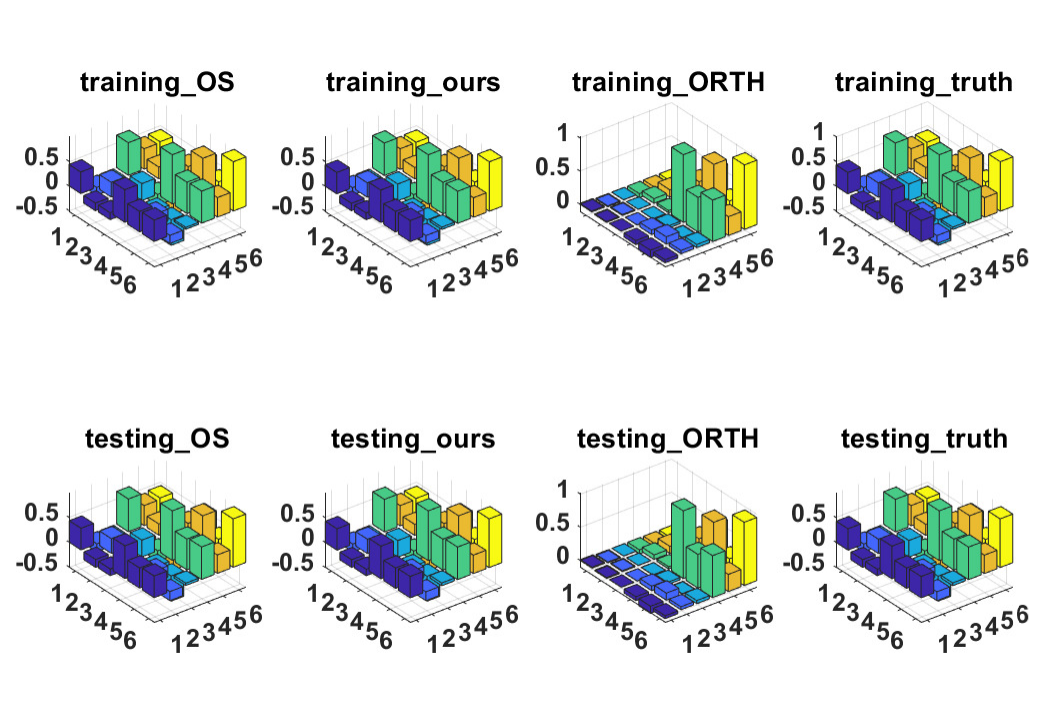}
    \caption{Measurement reconstruction covariance.}
    \label{fig:mea_rec}
\end{minipage}\hfill
\begin{minipage}{.5\textwidth}
    \centering
    \includegraphics[scale=0.40]{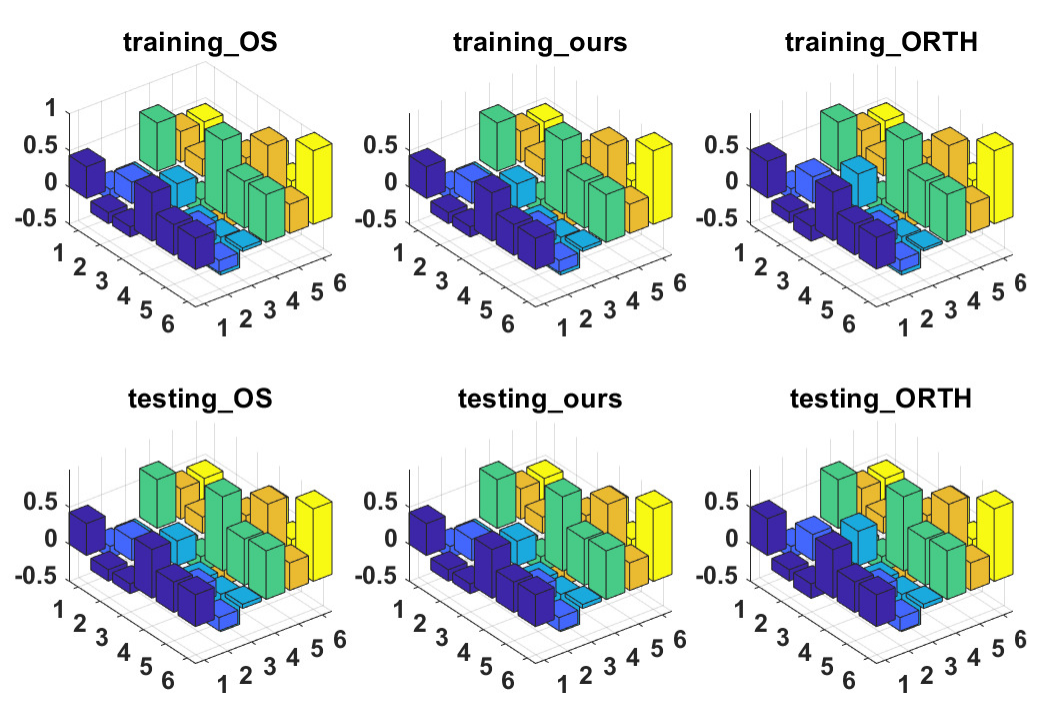}
    \caption{Measurement prediction covariance.}
    \label{fig:mea_pred}
\end{minipage}\\
\begin{minipage}{.5\textwidth}
    \centering
    \includegraphics[scale=0.40]{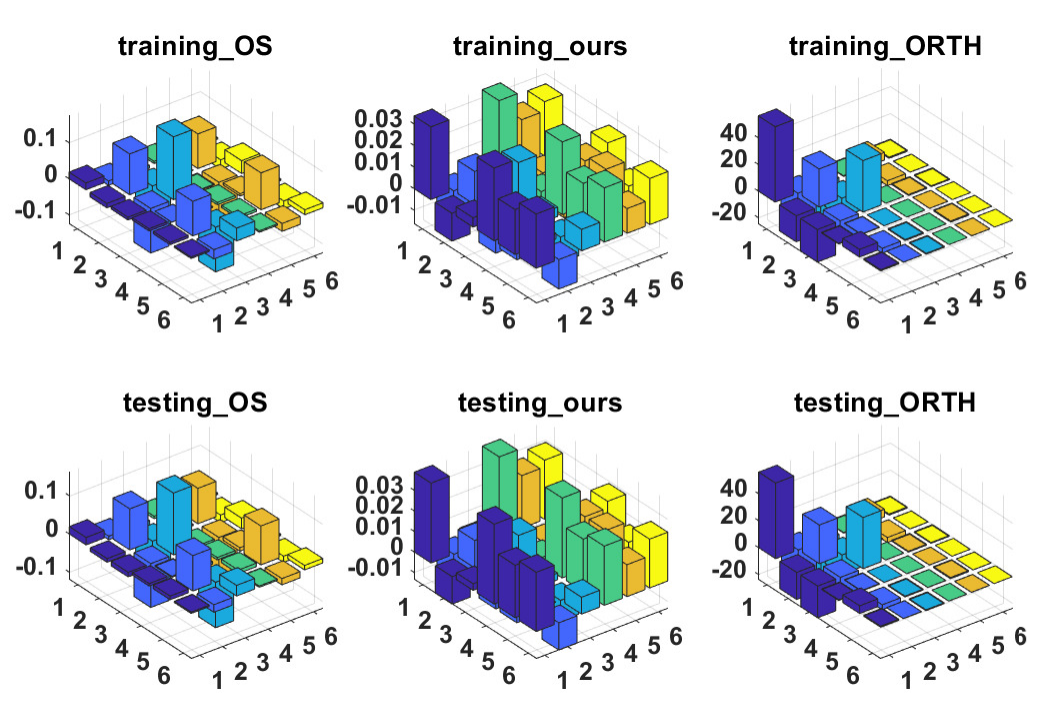}
    \caption{Signal reconstruction covariance. } 
    \label{fig:sig_rec}
\end{minipage}\hfill
\begin{minipage}{.5\textwidth}
    \centering
    \includegraphics[scale=0.40]{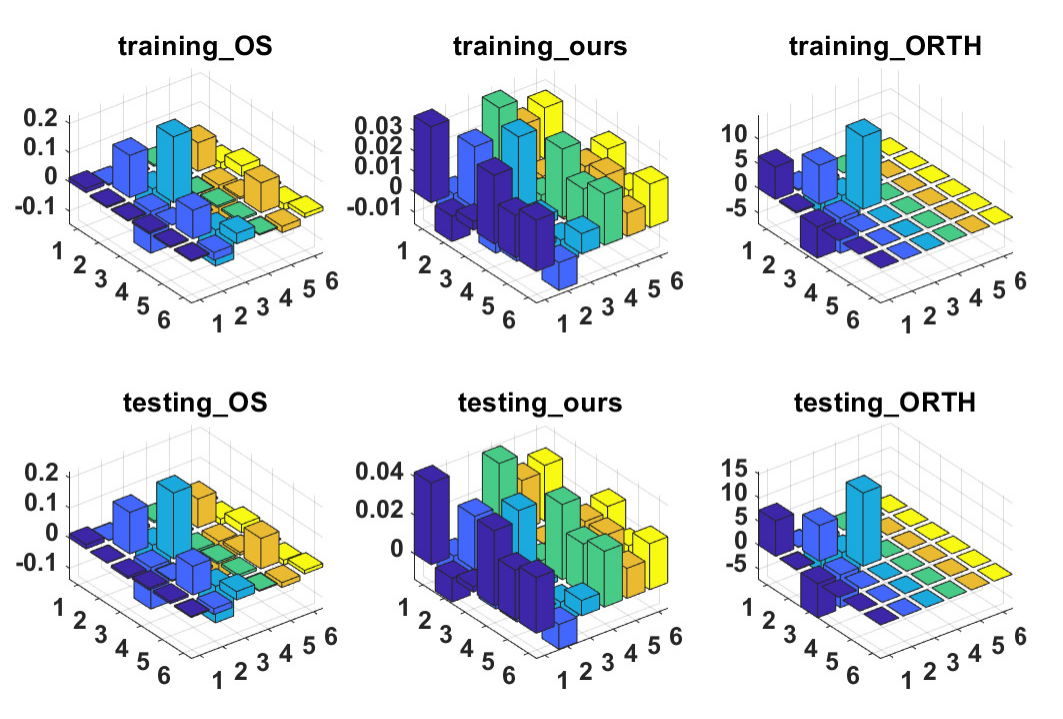}
    \caption{Signal prediction covariance.}
    \label{fig:sig_pred}
\end{minipage}
\end{figure*}

\begin{figure}[t]
    \centering
    \includegraphics[scale=0.40]{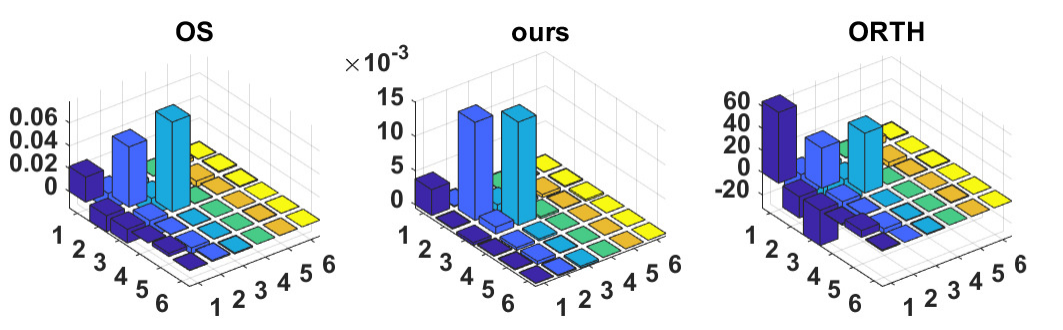}
    \caption{Signal prediction covariance by EM (non-orthogonal generation).}
    \label{fig:est_sig_pred}
\end{figure}

\subsection{Data Description}
The Lorenz oscillator described in~\cite{qin2022latent} is used to generate the latent dynamics of synthesized data. Note that the Lorenz oscillator with three coordinates~\mbox{$\bm v\in\Re^3$} is a nonlinear chaotic system instead of following the regression dynamics as in~\eqref{model_dy}. Assume that there are six sensors and each sensor measures a mix of the DLVs~$\bm v \in \Re^3$ and static noise~$\bm{\bar \varepsilon}\in\Re^3$ injected via respective channels.

Collect~$10000$ subsequent data points from the Lorenz oscillator and set the variance of the noise~$\bm{\bar \varepsilon}$ as that of the collected data. The loadings matrices are set as
\begin{align*}
    &\mathbf{P} = \begin{bmatrix}
    1&0&0\\
    0&1&0\\
    0&0&1\\
    0&0&0\\
    0&0&0\\
    0&0&0
    \end{bmatrix}\text{ and } \mathbf{\bar{P}} = \begin{bmatrix}
    $-$0.2997&$-$0.4611&$-$0.2868\\
    $-$0.2403&~0.2559&~0.6444\\
    $-$0.1334&~0.5749&$-$0.5168\\
    $-$0.2997&$-$0.4611&$-$0.2868\\
    $-$0.5400&$-$0.2052&~0.3576\\
    $-$0.6733&~0.3697&$-$0.1592
    \end{bmatrix}.
\end{align*}

The signal and noise subspaces are highly oblique as the canonical angles~\cite{qiu2005unitarily} between the subspaces spanned by the columns of~$\mathbf{R}$ and~$\mathbf{P}$ are~$23.99\degree$, $51.27\degree$, and~$60.97\degree$, which should be zeros for orthogonal cases. Obtain measurement samples~$\{\bm y_k\}$ by~\eqref{eq:outer_model}. Unless otherwise specified,  the first and last~$3000$ samples are for training and testing, respectively.

\begin{figure}
    \centering
    \includegraphics[scale=0.40]{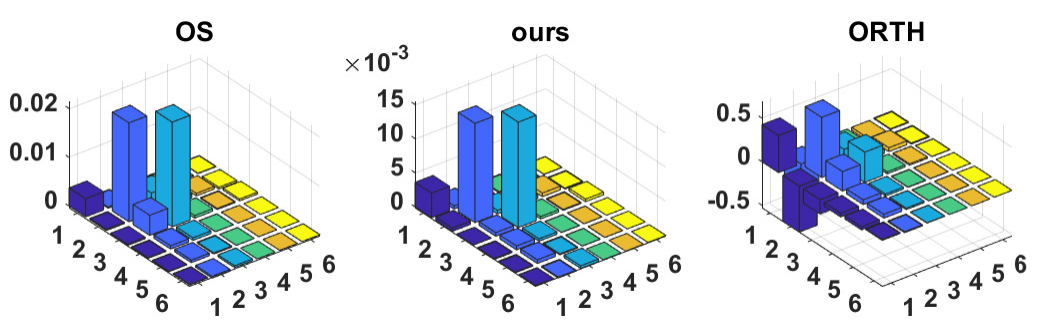}
    \caption{Signal prediction covariance by EM (ORTH generation).}
    \label{fig:orth_est_sig_pred}
\end{figure}
\subsection{Benchmark Algorithms}
The first benchmark (denoted as~OS) is called a one-shot algorithm. First, the eigendecomposition method in~\cite{gao2021modeling} is used to identify~$\mathbf{P}$ and~$\mathbf{R}$. Then, with the estimated oblique projection, update the parameters on the DLV dynamics as in Lines~\ref{alg: updateDLV} and~\ref{alg: updateB} of Algorithm~\ref{alg: DP}. As the OS algorithm does not have an alternate updating procedure, the comparison between the OS algorithm and ours will show the improving effect of the dynamics-projection interaction.

The second benchmark (denoted as~ORTH) emphasizes the orthogonal projection on the signal subspace. That is, the DLVs~$\bm v_k$ are estimated via the natural filter:
\begin{equation*} \label{or_filter}
   {\bm v}^{or}_k=(\mathbf{P}^\intercal \mathbf{P})^{-1}\mathbf{P}^{\intercal}\bm y_k = \mathbf{P}^{\dag}\bm y_k,
\end{equation*}
where~$\mathbf{P}^{\dag}$ is the Moore-Penrose inverse. The alternate updating procedure is still applicable and used. The comparison between this algorithm and ours will verify the improving effect of the oblique projection perspective.  


\begin{figure}[t]
    \centering
    \includegraphics[scale=0.21]
    {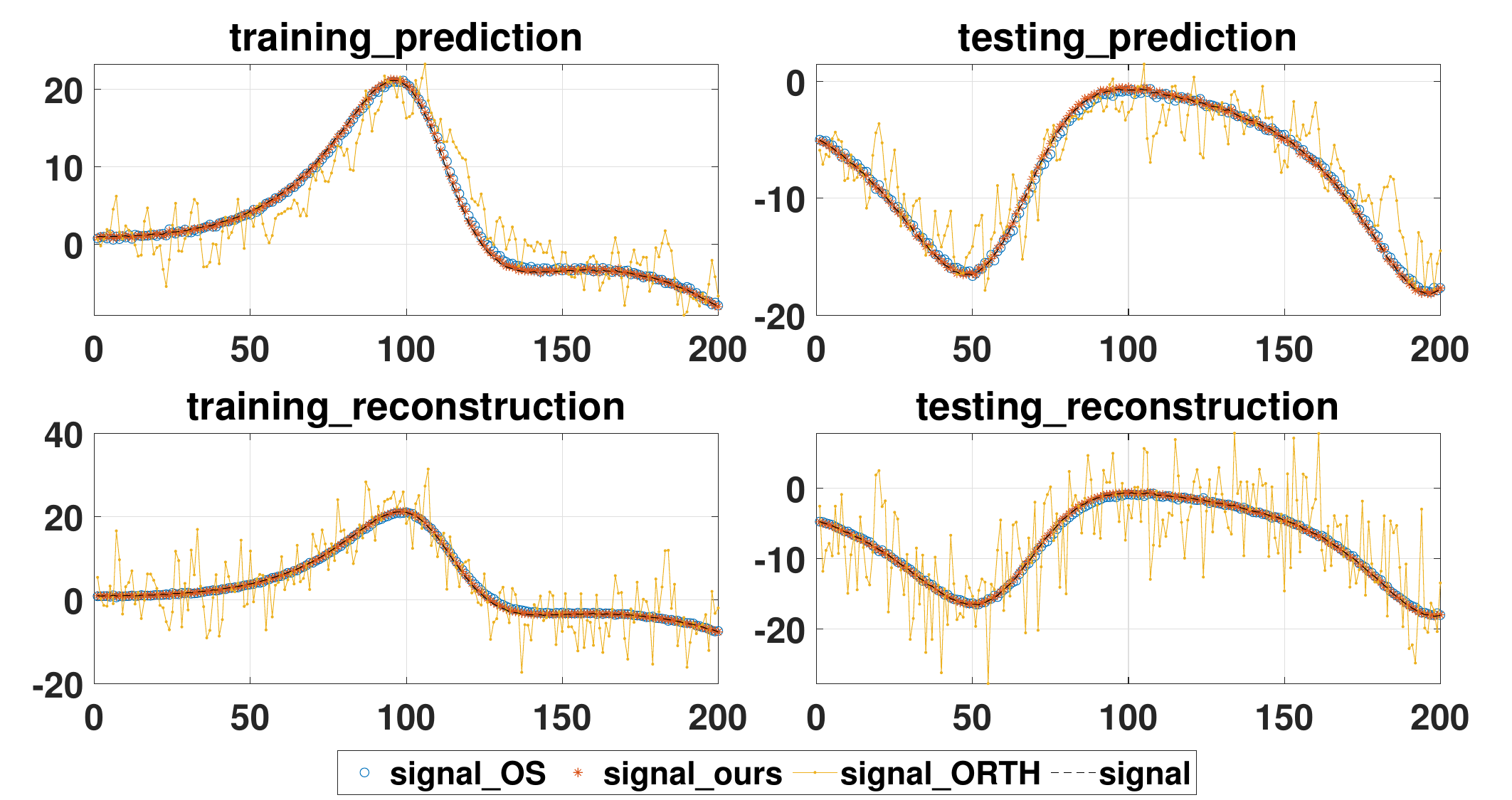}
    \includegraphics[scale=0.21]{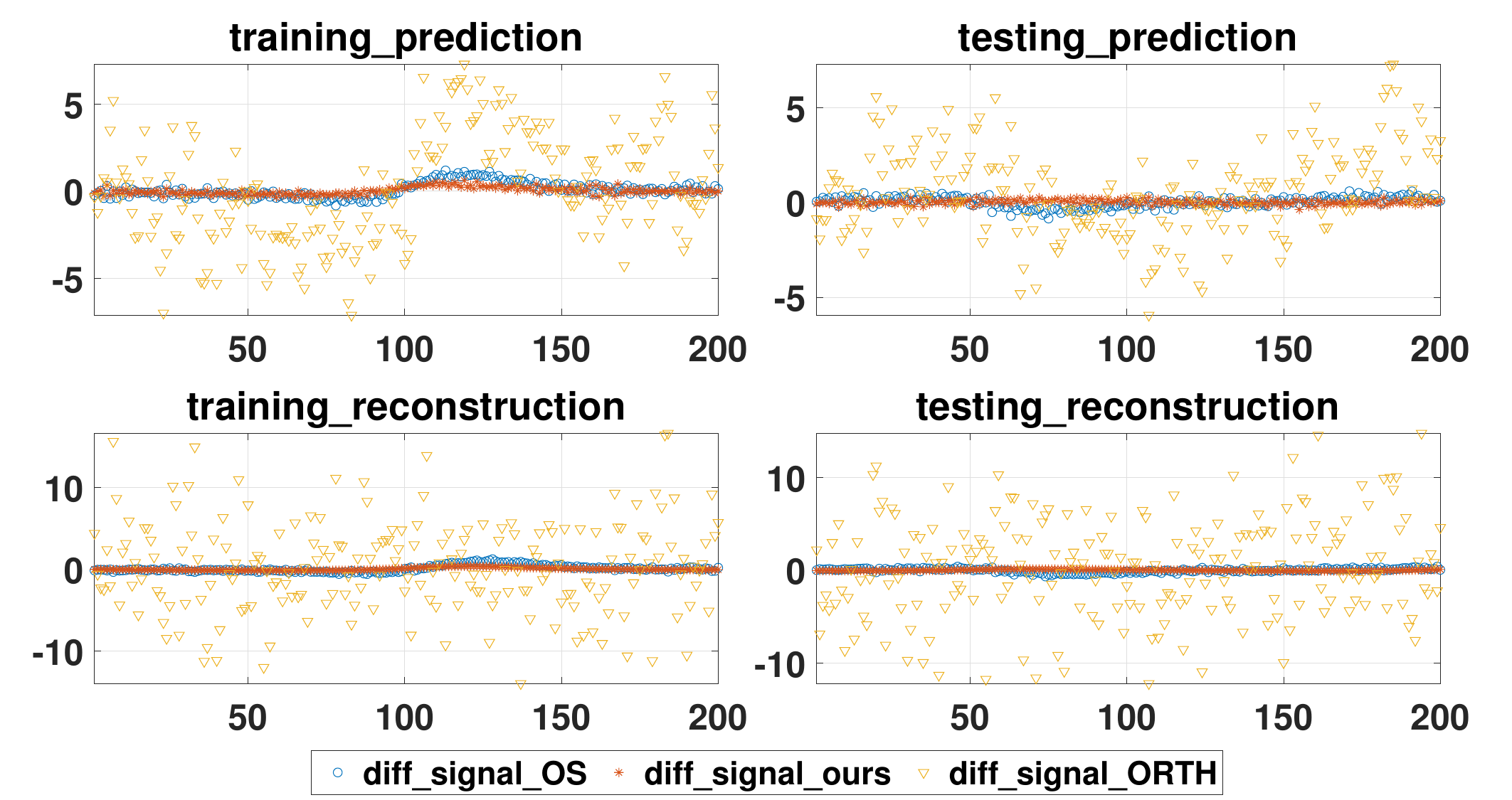}
    \caption{Reconstruction and prediction for Sensor 2.}
    \label{fig:sensor2}
\end{figure}

\subsection{Comparison}
The simulation results in the training and testing sets are presented separately. Similar performances regarding the two sets are observed for each algorithm. This fact suggests the efficacy of all algorithms. Nevertheless, the three algorithms have different efficiency levels.

Fig. \ref{fig:mea_rec} concerns the measurement reconstruction and depicts the covariance matrix of the series~$\{\bm y_k-\hat{\mathbf{P}}\hat{\mathbf{R}}^{\intercal}\bm y_k\}$. For the reference purpose, the covariance matrix (denoted as truth) of~$\{\bm y_k - \mathbf{P}\bm v_k\}$ is visualized, referring to the difference between the measurement and the true signal. As can be seen, the measurement reconstruction covariance for the OS algorithm or ours is similar to the reference covariance. However, the ORTH algorithm leads to a measurement reconstruction covariance matrix with smaller values, especially for those involving the first three sensors. This observation indicates that the ORTH algorithm does not eliminate sufficient static noise from the measurement for DLV dynamics identification, compared with two others. Fig. \ref{fig:mea_pred} depicts the covariance matrix of~$\{\bm y_k-\hat{\mathbf{P}}\hat{\bm{\tilde v}}_{k}\}$ and shows little difference among the three algorithms regarding measurement prediction. Combing this observation with~Fig. \ref{fig:mea_pred}, we argue that more noise is involved in the DLV dynamics estimation of the ORTH algorithm than the other two counterparts.

Fig. \ref{fig:sig_rec} depicts the covariance matrix of~$\hat{\mathbf{P}}\hat{\bm v}_k-\mathbf{P}\bm v_k$, referring to the difference between the reconstructed and true signals. Fig.~\ref{fig:sig_pred} depicts the covariance matrix of~$\hat{\mathbf{P}}\hat{\bm{\tilde v}}_{k}-\mathbf{P}\bm v_k$, referring to the difference between the predicted and true signals. As can be seen, the proposed algorithm attains the best reconstruction and prediction performances regarding covariance, while the ORTH algorithm achieves the worst. Fig.~\ref{fig:est_sig_pred} depicts~$\hat{\mathbf{P}}\hat{\mathbf{\Sigma}}_{\bm \varepsilon}\hat{\mathbf{P}}^{\intercal}$, where~$\hat{\mathbf{\Sigma}}_{\bm \varepsilon}$ is estimated from the EM procedure in each algorithm. Again, the proposed algorithm achieves the best performance in signal prediction, and the ORTH algorithm is the worst by likelihood analysis. These observations suggest that the proposed algorithm is the best to eliminate the static noise from the measurement and improve the DLV dynamics estimation.

The reconstructed~($\hat{\mathbf{P}}\hat{\mathbf{R}}^{\intercal}\bm y_k$), predicted~($\hat{\mathbf{P}}\hat{\bm{\tilde v}}_{k}$), and true~($\mathbf{P}\bm v_k$) signals in Sensor~$2$ are plotted in Fig.~\ref{fig:sensor2}. Also, the difference between the reconstructed and true signals and the difference between the predicted and true signals are plotted. The OS algorithm and ours can reconstruct or predict the signal well, but ours attains better performance. This fact exhibits the strength of our alternating iterative procedure. Also, we see many spikes in the signal curves generated by the ORTH algorithm. This observation again verifies that too much noise is reserved in the DLV dynamics estimation of the ORTH algorithm. Nevertheless, the curves generated by the ORTH algorithm can track the signal curves. This fact is consistent with the observation in Fig.~\ref{fig:sig_rec} and Fig.~\ref{fig:sig_pred} that the signal prediction performance is better than the signal reconstruction performance for the ORTH algorithm.

\begin{figure}
    \centering
    \includegraphics[scale=0.21]{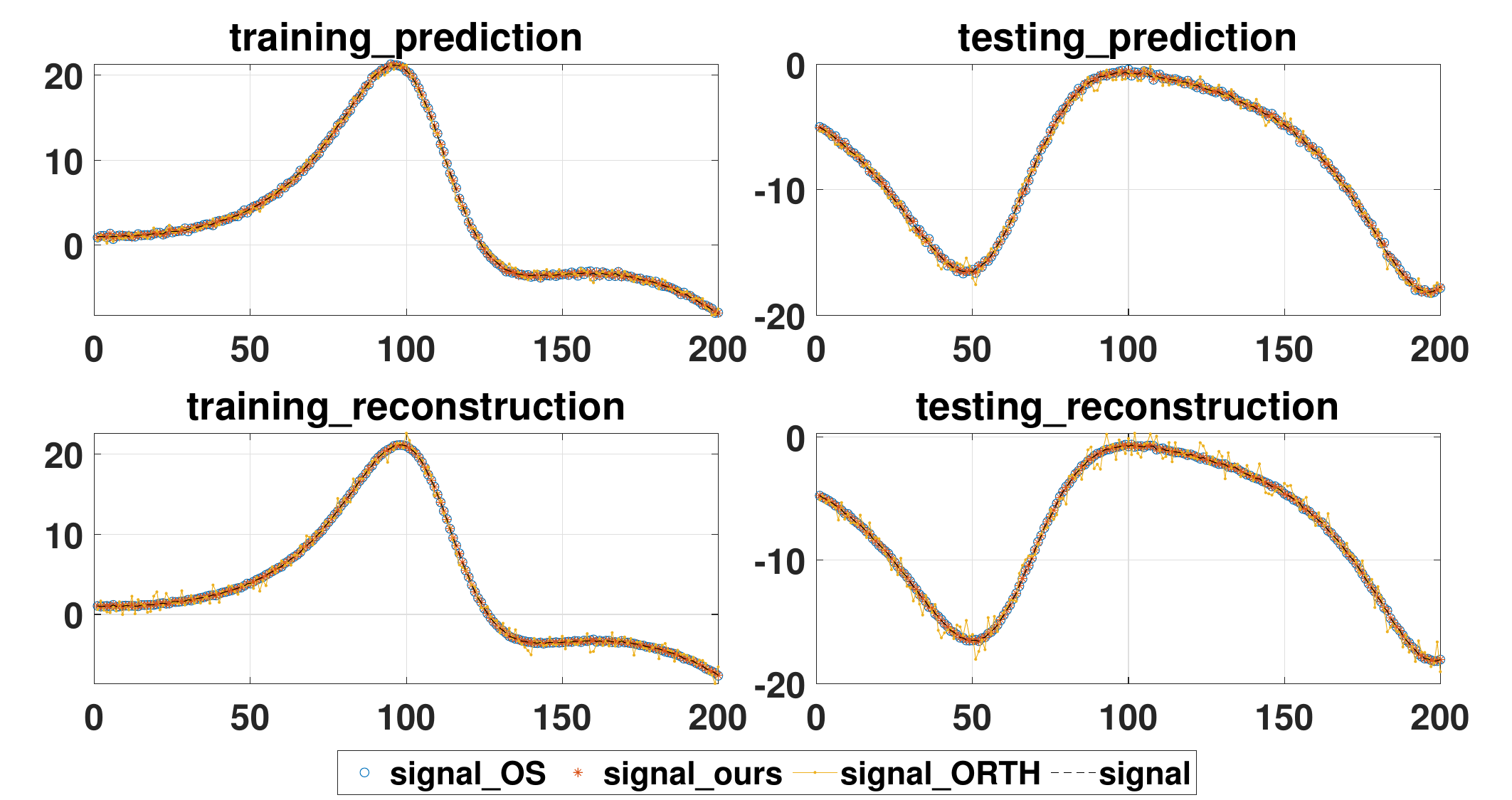}
    \includegraphics[scale=0.21]{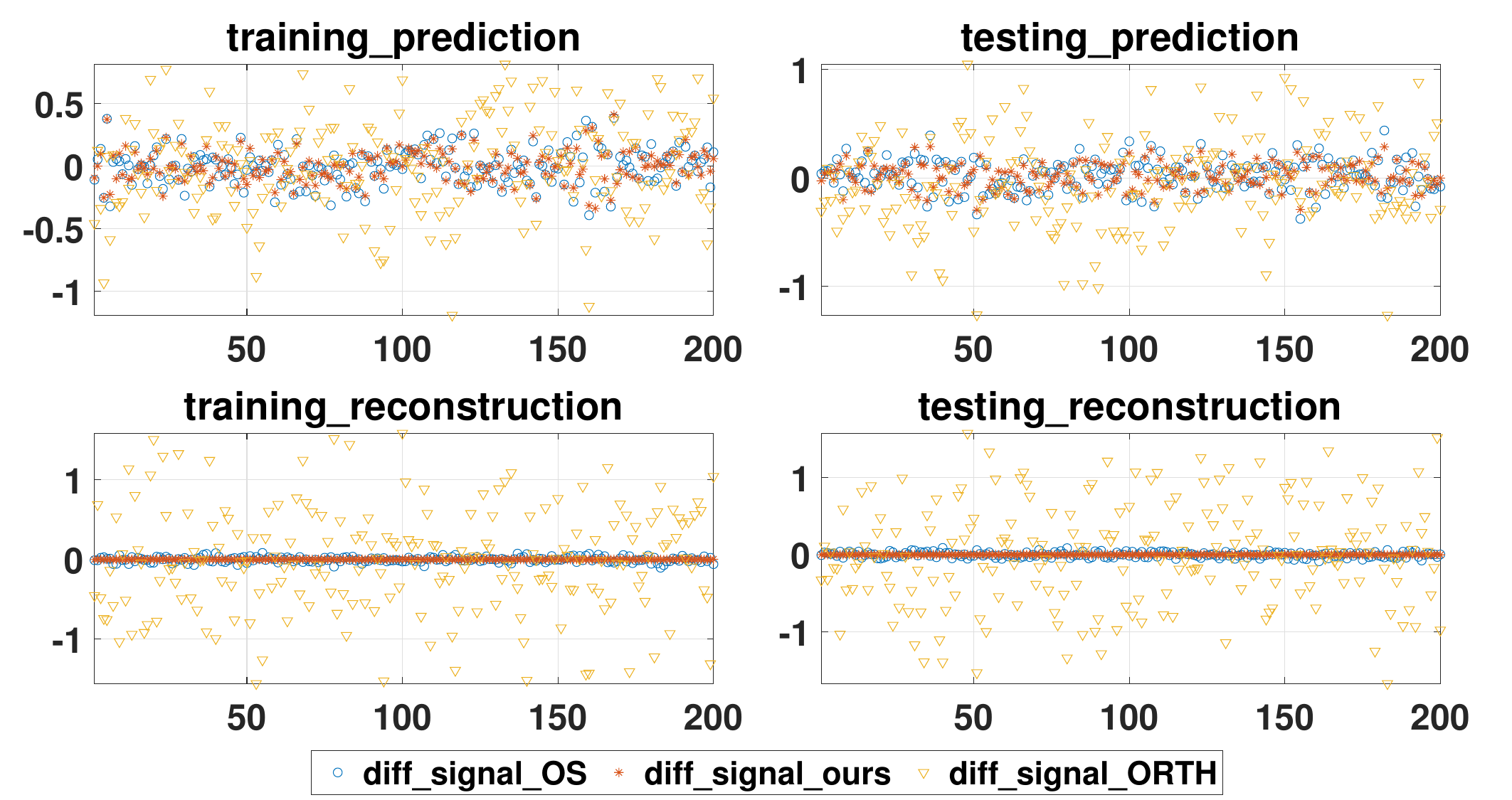}
    \caption{Reconstruction and prediction for Sensor 2 (ORTH generation).}
    \label{fig:orth_sensor2}
\end{figure}

\subsection{Empirical Consistency Analysis}
This subsection shows how the training sample set affects the oblique subspaces estimation. Table~\ref{tab: consistency} records the Frobenius norm of the difference between the true projection matrix~$\mathbf{P}\mathbf{R}^\intercal$ and its estimate from an algorithm by using the first~$x$ samples in the data set, with~$x$ varying from~1000 to 10000. Overall, the projection estimations of all three algorithms benefit from sample augmentation, though the effect diminishes and occasionally reverses as the sample number increases. Moreover, for the same training set, the proposed algorithm still attains a projection closest to the ground truth than the two benchmarks. We observe similar phenomena for signal subspace identification. Table~\ref{tab: consistency_signal} records the average canonical angle between the true signal subspace spanned by~$\mathbf{P}$ and the estimated signal subspace spanned by~$\hat{\mathbf{P}}$ for each algorithm. A smaller angle usually indicates a better signal subspace estimation. The two tables show the empirical consistency for all algorithms. Moreover, the proposed algorithm can extract more accurate information with fewer measurements, benefiting from the dynamics-projection interaction and the oblique projection perspective.

\begin{table}[t]
\centering
\caption{Distance to the true projection matrix}
\vspace{.5\baselineskip}
  \label{tab: consistency}
  \begin{tabularx}{0.68\linewidth}{c|ccc}
  \toprule    \toprule
 \diagbox{Num.}{Alg.} & OS &  \modelname & ORTH\\\hline
 1000& $1.2362$ & $0.8774$ & $2.0075$\\
 2000& $0.9113$ & $0.8367$ & $1.9381$\\
 3000& $0.9338$ & $0.8488$ & $1.7015$\\
 4000& $0.9016$ & $0.8352$ & $1.7338$\\
 5000& $0.9031$ & $0.8369$ & $1.7537$\\
 6000& $0.8679$ & $0.8270$ & $1.8011$\\
 7000& $0.8601$ & $0.8262$ & $1.7787$\\
 8000& $0.8589$ & $0.8313$ & $1.7774$\\
 9000& $0.8476$ & $0.8271$ & $1.7765$\\
 10000& $0.8571$ & $0.8330$ & $1.7610$\\
  \bottomrule \bottomrule
  \end{tabularx}
\end{table}
\begin{table}[t]
\centering
\caption{Angle-distance to the true signal subspace}
\vspace{.5\baselineskip}
  \label{tab: consistency_signal}
  \begin{tabularx}{0.7\linewidth}{c|ccc}
  \toprule    \toprule
 \diagbox{Num.}{Alg.} & OS & \modelname & ORTH\\\hline
 1000& $10.42\degree$ & $~4.31\degree$ & $22.68\degree$\\
 2000& $~5.30\degree$ & $~2.96\degree$ & $16.32\degree$\\
 3000& $~4.34\degree$ & $~3.05\degree$ & $~4.28\degree$\\
 4000& $~3.51\degree$ & $~1.76\degree$ & $~2.90\degree$\\
 5000& $~3.44\degree$ & $~1.49\degree$ & $~4.15\degree$\\
 6000& $~2.74\degree$ & $~1.19\degree$ & $~3.47\degree$\\
 7000& $~2.93\degree$ & $~1.38\degree$ & $~2.86\degree$\\
 8000& $~2.43\degree$ & $~1.32\degree$ & $~2.10\degree$\\
 9000& $~1.44\degree$ & $~0.94\degree$ & $~1.87\degree$\\
 10000& $~1.74\degree$ & $~1.02\degree$ & $~1.94\degree$\\
  \bottomrule \bottomrule
  \end{tabularx}
\end{table}

\subsection{Data with Orthogonal Signal and Noise Subspaces}
To further demonstrate the merit of the oblique-projection perspective, the measurements are generated using the same data points from the Lorenz oscillator and noise data~$\{\bm{\bar \varepsilon}_k\}$ but changing the static loadings matrix~$\mathbf{\bar{P}}$ as
\begin{align*}
    \begin{bmatrix}
        0&0&0&1&0&0\\
        0&0&0&0&1&0\\
        0&0&0&0&0&1
    \end{bmatrix}^{\intercal}.
\end{align*}
In this case, the signal and noise subspaces are orthogonal (ORTH generation). Fig.~\ref{fig:orth_est_sig_pred} depicts~$\hat{\mathbf{P}}\hat{\mathbf{\Sigma}}_{\bm \varepsilon}\hat{\mathbf{P}}^{\intercal}$, where~$\hat{\mathbf{\Sigma}}_{\bm \varepsilon}$ is estimated from the corresponding EM procedure in each algorithm. As seen, the proposed algorithm again pursues the smallest signal prediction covariance. Also, the reconstructed, predicted, true signals and their differences in Sensor~$2$ are plotted in Fig.~\ref{fig:orth_sensor2}. The performance of the ORTH algorithm improves in this case but is still worse than the other two. Notably, there are still spikes in the signal curves generated by the ORTH algorithm. The underlying reason may be that the Lorenz attractor is a nonlinear system, which leads to model mismatch with a linear VAR model to capture its dynamics. The ORTH algorithm is sensitive to the model mismatch. In contrast, the additional freedom can endow oblique projection with the ability to compensate for the DLV model mismatch against orthogonal projection that requires~$\mathbf{R}$ and~$\mathbf{P}$ to form the same subspace.



\bibliographystyle{ieeetr}
\bibliography{IEEEabrv,ProbLaVAR}
\end{document}